\documentclass[12pt,fleqn, a4]{article}

\usepackage[french]{babel}
\usepackage{inputenc}
\usepackage[T1]{fontenc}

\usepackage{amssymb}
\usepackage{latexsym}
\usepackage{graphics}

\usepackage{graphics}
\usepackage{color}
\newcommand{\colval}{0.3}
\definecolor{colone}{gray}{\colval}

\newcommand{\dcb}{\begin{array}{lll}}
\newcommand{\dce}{\end{array}}
\newcommand{\ebe}{\begin{enumerate}\setlength{\baselineskip}{13pt}\setlength{\parskip}{5pt}}
\newcommand{\dbe}{\end{enumerate}}

\newcommand{\ibegin}{\begin{itemize}\setlength{\baselineskip}{19pt}\setlength{\parskip}{7pt}}
\newcommand{\iend}{\end{itemize}}
\newcommand{\ok}{\rule{4pt}{6pt}}
\newcommand{\desb}{\begin{description}}
\newcommand{\dese}{\end{description}}

\newtheorem{Thm}{Theorem}[section]
\newtheorem {Cor}[Thm]{Corollary}
\newtheorem {definition}[Thm]{Definition}
\newtheorem {pro}{Proposition}[Thm]
\newtheorem {Lemma}[Thm]{Lemma}
\newtheorem {rem}[Thm]{Remark}

\newtheorem {assumption}[Thm]{Assumption}

\newcommand {\bd}{\begin{definition}}
\newcommand {\ed}{\end{definition}}
\newcommand {\bpro}{\begin{pro}}
\newcommand {\epro}{\end{pro}}
\newcommand {\bl}{\begin{Lemma}}
\newcommand {\el}{\end{Lemma}}
\newcommand {\bcor}{\begin{Cor}}
\newcommand {\ecor}{\end{Cor}}
\newcommand {\brem }{\begin{rem} \rm }
\newcommand {\erem }{\end{rem}}
\newcommand{\bethe}{\begin{Thm}}
\newcommand{\ethe}{\end{Thm}}
\newcommand {\bassumption}{\begin{assumption}}
\newcommand {\eassumption}{\end{assumption}}

\def \ind{1\!\!1}
\def\cro#1{\langle #1\rangle}

\begin{document}

\begin{center}
\textbf{\Large $\natural$-model with jumps}
\end{center}

\begin{center}
Shiqi Song

{\footnotesize Laboratoire Analyse et Probabilités\\
Université d'Evry Val D'Essonne, France\\
shiqi.song@univ-evry.fr}
\end{center}

\

\

\begin{quote}
\textbf{Abstract.} We consider the so-called $\natural$-model. It is an one-default model which gives the conditional law of a random time with respect to a reference filtration. This model has been studied in the case where the parameters are continuous. In this paper we will establish the $\natural$-model in the case of jump parameters. We then prove the corresponding enlargement of filtration formula and we compute the derivative of the conditional distribution functions of the random time.
\end{quote}

\section{Introduction}

We consider one-default model, i.e. the data of a random time $\tau$ combined with a filtration $\mathbb{F}$ under a probability measure $\mathbb{Q}$. The one-default models are widely applied in modeling financial risk and in price valuation of financial products such as CDS. The usefulness of an one-default model is conditional upon the way that the conditional laws of $\tau$ can be computed with respect to the filtration $\mathbb{F}$. The most used examples of random times, therefore, are the independent time, the Cox time, the honest time, the pseudo stopping time, the initial time, etc (cf., for example, \cite{lisbonn, BJR, EJY, JLC, JS1, LR, NP, NY}). In the paper \cite{JS2} a new class of random times has been introduced. Precisely, it is proved that, for any continuous increasing process $\Lambda$ null at the origin, for any continuous non-negative local martingale $N$ such that $0< N_te^{-\Lambda_t}< 1$, $t>0$, for any continuous local martingale $Y$, for any Lipschitz function $f$ on $\mathbb{R}$ null at the origin, there exist a random variable $\tau$ such that the family of conditional expectations $M^u_t=\mathbb{Q}[\tau\leq u|\mathcal{F}_t]$, $0< u,t<\infty,$ satisfy the following stochastic differential equation:
$$ 
(\natural_u)  
\left\{\dcb
dM^u_t=M^u_t\left(-\frac{e^{-\Lambda_t}}{1-N_te^{-\Lambda_t}}dN_t+f(M^u_t - (1-Z_t))dY_t\right),\ t\in[u,\infty),\\
M^u_u=1-N_ue^{-\Lambda_u}, \dce \right.
$$
We call this setting a $\natural$-model. 

There are two remarkable properties about the $\natural$-model. It is the only one in which the conditional laws of $\tau$ with respect to $\mathbb{F}$ are defined by a system of dynamic equations. The $\natural$-equation displays the evolution of the defaultable market. The knowledge of market evolution is a valuable property. This evolution form of the $\natural$-model had allowed \cite{JS2} to establish the so-called enlargement of filtration formula. It also is proved in \cite{JS2} that, reciprocally, the $\natural$-equation can be recovered from the enlargement of filtration formula in a way similar to that a differentiable function can be deduced from its derivative. 

We recall that the formula of enlargement of filtration is essential, when the no-arbitrage price valuation is considered in an one-default model (cf. \cite{ACDJ, FJS, songdrift}). Much as the enlargement of filtration formula is universally valid before the default time $\tau$, for a long time, the part of the enlargement of filtration formula after $\tau$ was merely proved for the honest time model or the initial time model. The $\natural$-models constitute the third family of models where the enlargement of filtration formula is valid on the whole $\mathbb{R}_+$. In addition, the enlargement of filtration formula in the $\natural$-model has a richer structure than that of honest time model, and has a more accurate expression than that of the initial time model. 

We recall also how widely the financial models are defined by stochastic differential equations, because it is one of the best ways to represent the evolution of a financial market. Usually, in an one-default model, there is no such a possibility to design the evolution. Now with the $\natural$-model, this becomes available. 

The second remarkable property of the $\natural$-model is its rich and flexible system of parameters $(Z,Y,f)$. The parameter $Z$ corresponds to the Azéma supermartingale of the $\natural$-model and determines the default intensity. The parameters $Y$ and $f$ describe the evolution of the market after the default time $\tau$. Such a system of parameters sets up a propitious framework for infering the market behavior and for calibrating the financial data.

We believe that the $\natural$-model can be a useful instrument to modeling financial market. We asked therefore if the $\natural$-model could be extended to the case of parameter processes with jumps. That is done. The purpose of this paper is to demonstrate the extension of the $\natural$-model to the jump case.

Here is an overview of the main results. 
\ebe
\item
The $\natural$-model is founded on the martingale point of view of the one-default models. It was introduced in \cite{JS2} and called i$\!M_Z$. In Section \ref{imz}, the results about i$\!M_Z$ are recalled with some modifications.
\item
The $\natural$-model with jumps will be defined in Section \ref{natural} for suitable parameters $Z,\mathbf{F},\mathbf{Y}$, where $Z$ is an Azema supermartingale, $\mathbf{F}$ is a Lipschitz functional in the sense of \cite{protter}, and $\mathbf{Y}$ is a local martingale. 

\item
One thing which can not be really explained in the text is how we have arrived at that definition. The first problem when we tried to extend the $\natural$-model from continuous case to jump case, was how the equation $(\natural)_u$ should be changed to take account of the jumps. After numerous essays, we conclude that the correct change is to replace the term $-\frac{e^{-\Lambda_t}dN_t}{1-N_te^{-\Lambda_t}}$ by $-\frac{dM_t}{{^\bullet}(1-Z_t)}$, where $Z$ is the given Azéma supermartingale; $M$ is the martingale part of $Z$; ${^\bullet}(1-Z_t)$ denotes the $\mathbb{F}$ predictable projection of $(1-Z)$. Our reflection is greatly influenced by \cite{MeYo, yoeurp} and by the fact that $-\frac{\Delta_t M}{{^{\mathbb{F}\cdot p}}(1-Z_t)}>-1$.

\item
The stochastic differential equation of the form $dX_t={^{\bullet}}X_tdW_t+dV_t$, where $W$ is a local martingale; $V$ is an increasing predictable process; ${^{\bullet}}X$ denotes the predictable projection of $X$, plays an important role in the study of the $\natural$-model. This equation had been considered in a different form in \cite{MeYo, yoeurp}. We obtain an explicite formula in Section \ref{natural} for the solution of such an equation, which is different than that of \cite{MeYo, yoeurp}. We wonder if this formula exists already in the literature. 

\item
To establish the $\natural$-model in the continuous case, we needed to look at the equation $(\natural)_u$ in its general form
$$ 
(\natural_u)  
\left\{\dcb
dX_t=X_t\left(-\frac{e^{-\Lambda_t}}{1-N_te^{-\Lambda_t}}dN_t+f(X_t - (1-Z_t))dY_t\right),\ t\in[u,\infty),\\
X_u=x.
\dce 
\right.
$$
We needed to prove three properties on the solutions $X^x$ of the equation $(\natural)_u$: $X^x_t\geq 0, X^x_t\leq 1-Z_t$ and $X^x_t\leq X^y_t$ if $0\leq x\leq y\leq 1-Z_u$. For the first property, it was the consequence of the Doleans-Dade exponential formula. The second property was proved by the local time technique (cf. \cite{RY}). Finally the third property is issued from the one-dimensional comparison theorem. It happens that the local time technique and the comparison theorem become inefficient in the jump case. This technical difficulty is overcomed by introducing the notion of $\natural$-pair in Section \ref{natural}.

\item
In Section \ref{efmc} we prove the enlargement of filtration formula for the $\natural$-models with Markovian coefficients. A striking point in this section is the use of the prediction process (cf. \cite{aldous, knight, MYr, yor}). The notion of the prediction process was introduced to represent the filtrations as Markov processes, which would give a pleasant way to make calculus on the filtrations. For long, the prediction process was assumed to contribute to the theory of enlargement of filtrations. See the study given in \cite{songthesis}. Nevertheless, the prediction process had not been widely applied in the literature, because likely the prediction process seemed not indispensable in the known examples. (For example, it was not used in the study of the continuous $\natural$-model in \cite{JS2}.) However, with the presence of the jumps, the prediction process appears unavoidable in the establishment of the enlargement of filtration formula in the $\natural$-models.
  
\item
When a one-default model is applied for a practice purpose, explicite formula for the conditional laws $\mathbb{Q}[\tau\in du|\mathcal{F}_t]$ will be needed. In Section \ref{conditionallaw} we explain how to compute these conditional laws in $\natural$-models. Actually in the case of Markovian coefficient with regularity, we will be able to compute the derivatives$$
\lim_{u\rightarrow} \frac{\mathbb{Q}[u<\tau\leq v|\mathcal{F}_t]}{A_v-A_u},
$$
where $A$ is the drift part of the submartingale $1-Z$. This yields in particular the absolute continuity of $\mathbb{Q}[\tau\in du|\mathcal{F}_t]$ for finite $t$ with respect to $dA_u$.

\item
One can wonder if the $\natural$-model is a particular case of known models such as honest model or initial models. The answer is no. However, many links exist between $\natural$-model and the other models. Roughly speaking, when $Z=e^{-\Lambda}$ for some continuous increasing predictable process $\Lambda$ with $\Lambda_0=0$ and when $\mathbf{F}\equiv 0$, the $\natural$-model will yield a Cox time. When the drift part of $1-Z$ is absolutely continuous with respect to Lebesgue measure and when $\mathbf{F}$ is of Markovian form, the $\natural$-model will yield an initial time. When $Z$ is a predictable decreasing process, the $\natural$-model will yield a pseudo stopping time. Certainly, if non of these conditions are satisfied, the $\natural$-model can produce original random times. Moreover, the $\natural$-model studied in this paper can not be honest time because of the assumption \textbf{Hy}$(Z)$ below.

\dbe

\

\section{Increasing family of bounded and positive martingales}\label{imz}

This section is borrowed from \cite[Section 2]{JS2} with some modifications.

\subsection{Product probability space}

As explained in \cite{JS2} an one-default model can always be imbeded isomorphically into a product probability space. Henceforth in this paper, models on product spaces only will be studied. Precisely, consider a measurable space $(\Omega,\mathcal{A})$ ($\mathcal{A}$ being a $\sigma$-algebra) equipped with a filtration $\mathbb{F}=(\mathcal{F}_t)_{t\geq 0}$ of sub-$\sigma$-algebras of $\mathcal{A}$, and consider the product space $[0,\infty]\times\Omega$ equipped with the product $\sigma$-algebra
$\mathcal{B}[0,\infty]\otimes\mathcal{A}$. Introduce the two maps $\pi,\tau$ defined as follows: $\pi(s,\omega)=\omega$ and $\tau(s,\omega)=s$ and introduce the filtration $\widehat{\mathbb{F}}=\pi^{-1}(\mathbb{F})$ on $[0,\infty]\times\Omega$. Provided with the pair $(\tau,\widehat{\mathbb{F}})$, constructing an one-default models on the product space amounts to constructing a probability measures on $\mathcal{B}[0,\infty]\otimes\mathcal{A}$. Recall that, working with the product space, it is custom to identify $\omega$ with the map $\pi$. In this way the filtration $\mathbb{F}$ is identified with $\widehat{\mathbb{F}}$ and the functions on $\Omega$ become functions on $[0,\infty]\times\Omega$. That is what we assume in this paper.

\subsection{i$\!M$ family associated with a probability measure $\mathbb{Q}$}

There exists a variety of ways to construct probability measures on the product space. Our approach in this paper is based on the following observation. A probability measure $\mathbb{Q}$ on $\mathcal{B}[0,\infty]\otimes\mathcal{A}$ is determined by its disintegration into its restriction on $\mathcal{A}$ and the conditional law of $\tau$ given $\mathcal{A}$. We will consider this conditional law $\mathbb{Q}[\tau\in du|\mathcal{A}]$ as the terminal value of the measure-valued martingale $(\mathbb{Q}[\tau\in
du|\mathcal{F}_t]:0\leq t\leq \infty)$ and we will define this measure-valued martingale by a stochastic differential equation. For this purpose, we introduce the following notion.

Let $\mathbb{P}$ be a probability measure on $\mathcal{A}$. An increasing family of positive $(\mathbb{P},\mathbb{F})$ martingales bounded by
$1$ (in short i$\!M(\mathbb{P},\mathbb{F})$ or simply i$\!M$) is a family of processes $(M^u: 0\leq
u\leq \infty)$ satisfying the following conditions: 
\ebe 
\item 
Every $M^u$ is a càdlàg non-negative $(\mathbb{P},\mathbb{F})$ martingale on
$[u,\infty]$, bounded by $ 1$, and closed by $M^u_\infty$. 
\item 
For every $0\leq t\leq \infty$, the random map $u\in[0,t]\rightarrow M^u_t$ is a right continuous non-decreasing function. 
\item 
$M^\infty_\infty=1$.
\dbe

The theorem below gives the relationship between an i$\!M$ and the construction of a probability measure on the product space. We recall that we identify the elements on $(\Omega,\mathcal{A})$ with elements on the product space.

\bethe \label{theoremIM} 
Let $\mathbb{P}$ be a probability measure on $\mathcal{A}$. Suppose that the filtration $\mathbb{F}$ is right continuous and contains the $(\mathbb{P},\mathcal{F}_\infty)$ null sets. 
\ebe 
\item For any probability measure $\mathbb{Q}$ on the $\sigma$-algebra $\mathcal{B}[0,\infty]\otimes\mathcal{A}$ which coincides with $\mathbb{P}$ on $\mathcal{F}_\infty$, there exists a unique i$M(\mathbb{P},\mathbb{F})=(M^u:0\leq u\leq \infty)$ such that, for $0\leq u\leq t\leq  \infty$, $$
M^u_t=\mathbb{Q}[\tau\leq u|\mathcal{F}_t].
$$
We shall say that this i$M$ is associated with $\mathbb{Q}$.

\item 
Let $(M^u: 0\leq u\leq \infty)$ be an i$M(\mathbb{P},\mathbb{F})$. There is a
unique probability measure $\mathbb{Q}$ on
the $\sigma$-algebra $\mathcal{B}[0,\infty]\otimes\mathcal{A}$ which coincides with $\mathbb{P}$ on $\mathcal{A}$ and satisfies $\mathbb{Q}[\tau\leq u|\mathcal{F}_t]=M^u_t$ for $0\leq
u\leq t\leq  \infty$, and $\mathbb{Q}[\tau\leq u|\mathcal{A}]=M^u_\infty$. We shall say that $\mathbb{Q}$ is
associated with the i$M$ and with $\mathbb{P}$.
\dbe
\ethe
\textbf{Proof.} Consider the first assertion. The uniqueness is clear. For each $0\leq
u\leq \infty$, let $(G^u_t: 0\leq t\leq \infty)$ be a càdlàg version of the
$(\mathbb{P},\mathbb{F})$ martingale $\mathbb{Q}[\tau\leq
u|\mathcal{F}_t], t\in[0,\infty]$. For $u<v$, for any $\mathbb{F}$ stopping time $T$, $\mathbb{P}$
almost surely, we have $0\leq G^u_T\leq G^v_T\leq 1$. Set $$ M^u_t=\inf_{w\in\mathtt{Q}_+,
w>u}(G^w_t\wedge 1)^+,\ 0\leq u<\infty,0\leq t\leq \infty.
$$
Then,
\begin{itemize}
\item 
For $0\leq u<\infty$, the process $M^u=(M^u_t)_{0\leq t\leq \infty}$ is $\mathbb{F}$-optional.
\item  
For $0\leq t\leq \infty$, $u\in[0,\infty)\rightarrow M^u_t$ is right continuous.
\item  
For $0\leq t\leq \infty, 0\leq u<v<\infty$, $0\leq M^u_t\leq M^v_t\leq 1$ everywhere.
\item 
For $0\leq u<\infty$, for any $\mathbb{F}$ stopping time $T$, $\mathbb{P}$
almost surely, $G^u_T\leq M^u_T$.
\end{itemize}
Let $0\leq u<\infty$ and $T$ be an $\mathbb{F}$ stopping time. We write$$ 
\dcb 
\mathbb{Q}[\tau\leq
u]=\mathbb{Q}[G^u_T]\leq\mathbb{Q}[M^u_T] =
\inf_{v\in\mathtt{Q}_+, u<v}\mathbb{Q}[G^v_T]
=\inf_{v\in\mathtt{Q}_+, u<v}\mathbb{Q}[\tau\leq v]
=\mathbb{Q}[\tau\leq u].
\dce
$$This shows that $M^u_T=G^u_T$, $\mathbb{P}$ almost surely. Consequently $G^u$ and $M^u$ are $\mathbb{P}$ indistinguishable (cf. \cite{Yan}). Define finally $M^\infty_t=1$ for $0\leq t\leq \infty$. The family of processes $(M^u:0\leq u\leq \infty)$ thus defined satisfies the conditions of an i$\!M$.

Consider the second assertion. The uniqueness is clear. With the càdlàg non-decreasing map $u\in[0,\infty]\rightarrow M^u_\infty$ we denote by $d_uM^u_\infty$ the associated random measure on $[0,\infty]$ (noting that at $0$ it can have a mass equal to $M^0_\infty$). Define a
probability measure on
$([0,\infty]\times\Omega,\mathcal{B}\otimes\mathcal{A})$ by
$$
\mathbb{Q}[F]=\mathbb{P}[\int_{[0,\infty]}F(t,\cdot)d_tM^t_\infty]
$$
where $F(t,\omega)\in\mathcal{B}[0,\infty]\otimes\mathcal{A},\ F(t,\omega)\geq 0$. We compute. First for $A\in\mathcal{A}$: $$
\mathbb{Q}[A]=\mathbb{Q}[A \cap \{0\leq \tau\leq \infty\}]
=\mathbb{P}[\mathbb{I}_A\int_{[0,\infty]}d_tM^t_\infty]
=\mathbb{P}[A].
$$Secondly let  $0\leq u\leq t\leq \infty, A\in\mathcal{F}_t$. We have $$
\mathbb{Q}[A \cap \{\tau\leq
u\}]=\mathbb{P}[\mathbb{I}_A\int_{[0,u]}d_sM^s_\infty]
=\mathbb{P}[\mathbb{I}_AM^u_\infty]=\mathbb{P}[\mathbb{I}_AM^u_t]
=\mathbb{Q}[\mathbb{I}_AM^u_t].
$$
This implies that $\mathbb{Q}[\tau\leq
u|\mathcal{F}_t]=M^u_t$. If $A\in \mathcal{A}$, we obtain $\mathbb{Q}[\tau\leq
u|\mathcal{A}]=M^u_\infty$.\ok

\subsection{i$\!M_Z$ family}

In this subsection, besides the given probability structure $(\mathbb{P},\mathbb{F})$, we consider a $(\mathbb{P},\mathbb{F})$ supermartingale $Z$ such that $0\leq Z\leq 1$. We introduce the following definition:

An increasing family of positive martingales bounded by $1-Z$ (in short i$\!M_Z(\mathbb{P},\mathbb{F})$ or simply i$\!M_Z$) is an i$\!M=(M^u: 0\leq u\leq\infty)$ satisfying the following conditions: for any $0\leq u\leq t<\infty$, $M^u_u=1-Z_u$ and $M^u_t \leq  1-Z_t$.

The theorem below is an immediate consequence of the Theorem \ref{theoremIM}.

\bethe \label{martingaleM} 
Let $\mathbb{P}$ be a probability measure on $\mathcal{A}$. Suppose that the filtration $\mathbb{F}$ is right continuous and contains the $(\mathbb{P},\mathcal{F}_\infty)$ null sets. Let $(M^u: 0\leq u\leq\infty)$ be an i$M(\mathbb{P},\mathbb{F})$ associated with a probability measure $\mathbb{Q}$ on $\mathcal{B}[0,\infty]\otimes\mathcal{A}$ which coincides with $\mathbb{P}$ on $\mathcal{F}_\infty$. Then, i$\!M$ is a i$\!M_Z$ if and only if $\mathbb{Q}[t<\tau|\mathcal{F}_t]=Z_t$ for $t\geq 0$.
\ethe

\

\section{$\natural$-model}\label{natural}

From now on we fix a stochastic structure $(\mathbb{P},\mathbb{F})$ satisfying the usual condition and a $(\mathbb{P},\mathbb{F})$ supermartingale $Z$ satisfying the following assumption.

\textbf{Hy}($Z$): For $0<t<\infty$, $1-Z_t>0$ and $1-Z_{t-}>0$.

\

\subsection{An affine stochastic differential equation for positive submartingales}

\bl\label{separable}
Let $u\geq 0$. Let $W$ be a $(\mathbb{P},\mathbb{F})$ local martingale with $\Delta W> -1$, $V$ be a non-decreasing $\mathbb{F}$ predictable process and $a\geq 0$. Consider the stochastic differential equation$$
\left\{
\dcb
d\Delta_t&=&{^{\mathbb{F}\cdot p}}\Delta_{t}dW_t+dV_t,\ t\in [u,\infty),\\
\Delta_u&=&a,
\dce
\right.
$$ 
where the superscript\ \ ${^{\mathbb{F}\cdot p}}$ denotes the predictable projection in $\mathbb{F}$. Then, the solution $\Delta^a$ of this equation is given by$$
\Delta^a_t=\mathcal{E}(\ind_{(u,\infty)}\centerdot W)_t\left(a+\int_u^t\frac{1}{\mathcal{E}(\ind_{(u,\infty)}\centerdot W)_{s-}}dV_s\right),\ t\in [u,\infty).
$$
\el

\textbf{Proof.} The stochastic differential equation in this lemma has uniqueness of the solution. Let $X$ to be the right hand term in the above formula. We apply the integration by parts formula to check that $X$ is the solution of the equation:
$$
\dcb
dX_t
&=&\mathcal{E}(\ind_{(u,\infty)}\centerdot W)_{t-}\left(a+\int_{(u,t)}\frac{1}{\mathcal{E}(\ind_{(u,\infty)}\centerdot W)_{s-}}dV_s\right) \ind_{(u,\infty)}dW_t + \ind_{(u,\infty)}dV_t\\
&&\hspace{4cm}+d\left[\mathcal{E}(\ind_{(u,\infty)}\centerdot W),\ind_{(u,\infty)}\frac{1}{\mathcal{E}(\ind_{(u,\infty)}\centerdot W)_{-}}\centerdot V\right]_t,\\

&=&\mathcal{E}(\ind_{(u,\infty)}\centerdot W)_{t-}\left(a+\int_{(u,t)}\frac{1}{\mathcal{E}(\ind_{(u,\infty)}\centerdot W)_{s-}}dV_s\right) \ind_{(u,\infty)}dW_t + \ind_{(u,\infty)}dV_t\\
&&\hspace{4cm}+\mathcal{E}(\ind_{(u,\infty)}\centerdot W)_{t-}\frac{1}{\mathcal{E}(\ind_{(u,\infty)}\centerdot W)_{t-}}\ind_{(u,\infty)}d\left[W,V\right]_t,\\

&=&\mathcal{E}(\ind_{(u,\infty)}\centerdot W)_{t-}\left(a+\int_{(u,t)}\frac{1}{\mathcal{E}(\ind_{(u,\infty)}\centerdot W)_{s-}}dV_s\right) \ind_{(u,\infty)}dW_t + \ind_{(u,\infty)}dV_t\\
&&\hspace{4cm}+\mathcal{E}(\ind_{(u,\infty)}\centerdot W)_{t-}\frac{1}{\mathcal{E}(\ind_{(u,\infty)}\centerdot W)_{t-}}\ind_{(u,\infty)}\Delta_tV dW_t,\\

&=&\mathcal{E}(\ind_{(u,\infty)}\centerdot W)_{t-}\left(a+\int_u^t\frac{1}{\mathcal{E}(\ind_{(u,\infty)}\centerdot W)_{s-}}dV_s\right) \ind_{(u,\infty)}dW_t + \ind_{(u,\infty)}dV_t\\
\\
&=&{^{\mathbb{F}\cdot p}}X_t\ind_{(u,\infty)}dW_t + \ind_{(u,\infty)}dV_t.\ \ok
\dce
$$

As a corollary we have:

\bl\label{positive}
Suppose the same setting as in the previous lemma. Then, the solution $\Delta^a$ is non-negative for all $t\in[u,\infty)$. If a>0, $\Delta^a>0$ on $[u,\infty)$. Let $T=\inf\{s\geq u: V_s-V_u>0\}$. Then $\Delta^0>0$ on $(T,\infty)$.
\el

We have a partially inverse result, which is a direct consequence of \cite[Théorème(6.31)]{jacod} and of \cite{MeYo}.

\bl\label{prv}
Let $u\geq 0$. Let $X$ be a $(\mathbb{P},\mathbb{F})$ submartingale such that $X>0$ and $\frac{1}{{^{\mathbb{F}\cdot p}}X}$ is locally bounded on $[u,\infty)$. Then there exists a $(\mathbb{P},\mathbb{F})$ local martingale $W$ with $\Delta W>-1$ and a non-decreasing $\mathbb{F}$ predictable process $V$ such that$$
dX_t={^{\mathbb{F}\cdot p}}X_{t}dW_t+dV_t,\ t\in[u,\infty).
$$
\el

\

\subsection{The positive submartingale $1-Z$}

We recall \textbf{Hy}(Z). Let $Z=M-A$ be the $(\mathbb{P},\mathbb{F})$ canonical decomposition of $Z$ with $M$ a $(\mathbb{P},\mathbb{F})$ local martingale and $A$ a non-decreasing $\mathbb{F}$ predictable process. Then, ${^{\mathbb{F}\cdot p}}(1-Z)_t=1-Z_{t-}+\Delta_tA>0$ for any $0<t<\infty$. For any $0<u<\infty$, the submartingale $1-Z$ satisfies the stochastic differential equation:$$
d(1-Z)_t
= {^{\mathbb{F}\cdot p}}(1-Z)_t\ \frac{-dM_t}{{^{\mathbb{F}\cdot p}}(1-Z)_t}+dA_t, \ t\in [u,\infty).
$$
Note that $-\Delta_t M+{^{\mathbb{F}\cdot p}}(1-Z)_t=1-Z_t$ so that $$
\frac{-\Delta_t M}{{^{\mathbb{F}\cdot p}}(1-Z)_t}=
\frac{1-Z_t}{{^{\mathbb{F}\cdot p}}(1-Z_t)}-1>-1, \ t\in [u,\infty).
$$
We define, for $0<u<\infty$, $\widetilde{m}^u_t=\int_u^t \frac{-dM_s}{{^{\mathbb{F}\cdot p}}(1-Z)_s}$. Since obviously $d\widetilde{m}^u_t=d\widetilde{m}^v_t$ for $0<u<v\leq t<\infty$, we omit the superscripts and we denote simply
$$
d\widetilde{m}_t=\frac{-dM_t}{{^{\mathbb{F}\cdot p}}(1-Z)_t},\ t\in(0,\infty).
$$
Below we denote simply $
{^{\mathbb{F}\cdot p}}X={^\bullet}X
$
for any process whose $\mathbb{F}$ predictable projection exists.

\

\subsection{$\natural$-equation and $\natural$-pair}\label{generatingEquation}

Let $\mathbb{D}$ design the space of all càdlàg $\mathbb{F}$ adapted processes. Let $m>0$ be an integer. Let $\mathbf{Y}=(Y_1,\ldots,Y_m)$ be an $m$-dimensional $(\mathbb{P},\mathbb{F})$ local martingale, and $\mathbf{F}=(F_1,\ldots,F_m)$ be a Lipschitz functional from $\mathbb{D}$ into the set of $m$-dimensional locally bounded $\mathbb{F}$ predictable processes in the sense of \cite{protter}. For $0< u<\infty$, for any ${\mathcal {F}}_u$-measurable random variable $x$, we consider the stochastic differential equation determined by the pair $(\mathbf{F},\mathbf{Y})$: $$
(\natural_u)  \left\{\dcb
dX_t&=&X_{t-} d\widetilde{m}_t+\mathbf{F}(X)_{t}^\top d\mathbf{Y}_t,\ t\in[u,\infty),\\
X_u&=&x.
\dce
\right.
$$

We will call the pair $(\mathbf{F},\mathbf{Y})$ a $\natural$-pair if it satisfies the following conditions, for any $1\leq j\leq m$, for any $u>0$ and for any $X,X'\in\mathbb{D}$:
\ebe
\item[(i)]
The process $t\in[u,\infty)\rightarrow\frac{F_j(X)_t}{{^\bullet}(1-Z)_t-X_{t-}}\ind_{\{(1-Z_{t-})-X_{t-}\neq 0\}}$ is integrable with respect to $Y_j$, and satisfies the inequality:
$$
\Delta_t\tilde{m}-\frac{1}{{^\bullet}(1-Z)_t-X_{t-}}\ind_{\{{^\bullet}(1-Z)_t-X_{t-}\neq 0\}}\mathbf{F}(X)_t^\top\Delta_t\mathbf{Y}>-1, \  t\in[u,\infty).
$$

\item[(ii)]
The process $t\in[u,\infty)\rightarrow\frac{F_j(X)_t}{X_{t-}}\ind_{\{X_{t-}\neq 0\}}$ is integrable with respect to $Y_j$, and satisfies the inequality:
$$
\Delta_t\tilde{m}+\frac{1}{X_{t-}}\ind_{\{X_{t-}\neq 0\}}\mathbf{F}(X)_t^\top\Delta_t\mathbf{Y}\geq -1, \  t\in[u,\infty).
$$

\item[(iii)]
If $0\leq X,X'\leq 1$, the process $t\in[u,\infty)\rightarrow\frac{F_j(X)_t-F_j(X')_t}{X_{t-}-X'_{t-}}\ind_{\{X_{t-}-X'_{t-}\neq 0\}}$ is integrable with respect to $Y_j$, and satisfies the inequality:
$$
\Delta_t\tilde{m}+\frac{1}{X_{t-}-X'_{t-}}\ind_{\{X_{t-}-X'_{t-}\neq 0\}}(\mathbf{F}(X)_t-\mathbf{F}(X')_t)^\top\Delta_t\mathbf{Y}\geq -1, \  t\in[u,\infty).
$$
\dbe

\brem
Note that the inequality in the condition (iii) is equivalent to$$
X_{t-}+X_{t-}\Delta_t\tilde{m}+\mathbf{F}(X)_t^\top\Delta_t\mathbf{Y}
\geq 
X'_{t-}+X'_{t-}\Delta_t\tilde{m}+\mathbf{F}(X')_t^\top\Delta_t\mathbf{Y},
$$
whenever $X_{t-}>X'_{t-}$. This condition looks like very much that one in \cite[Chapter V Section 10 Theorem 65]{protter} which ensures that the solutions of a stochastic differential equation form a diffeomorphism. 
\erem

\subsection{The basic properties of the $\natural$-equation }

Consider an $m$-dimensional $(\mathbb{P},\mathbb{F})$ local martingale $\mathbf{Y}=(Y_1,\ldots,Y_m)$ and a Lipschitz functional $\mathbf{F}=(F_1,\ldots,F_m)$ from $\mathbb{D}$ into the set of $m$-dimensional locally bounded $\mathbb{F}$ predictable processes in the sense of \cite{protter}.

\bl\label{three}
For $0<u<\infty$ let $X^x$ be the solution of the equation $(\natural_u)$ associated with the pair $(\mathbf{F},\mathbf{Y})$.  
\ebe 
\item[(1)]
Suppose that $X^x_u=x\leq 1-Z_u$. Then, if $(\mathbf{F},\mathbf{Y})$ satisfies the above condition (i) on $[u,\infty)$, we have $(1-Z)_t-X^x_t\geq 0$ for $t\in[u,\infty)$ and $(1-Z)_t-X^x_t> 0$ for $t\in(T,\infty)$ where $T=\inf\{s\geq u: (1-Z_u)-x+A_s-A_u> 0\}$. Inversely, for any $\mathbb{F}$ stopping time $u\leq S\leq T$, if $(1-Z)_{t}-X^x_{t}>0$ and $(1-Z)_{t-}-X^x_{t-}>0$ on $(S,T]$, the above condition (i) is satisfied on $(S,T]$ (instead of $[u,\infty)$) by $(\mathbf{F},\mathbf{Y})$ with $X=X^x$. 
\item[(2)]
Suppose that $X^x_u=x\geq 0$. Then, $X^x\geq 0$ and $\frac{1}{X^x_{-}}\ind_{\{X^x_->0\}}$ is $X^x$ integrable on $[u,\infty)$, if and only if the above condition (ii) is satisfied by $(\mathbf{F},\mathbf{Y})$ with $X=X^x$. In addition, $X^x>0, X^x_->0$ on $[u,\infty)$ if and only if $(\mathbf{F},\mathbf{Y})$ and $X=X^x$ satisfies the condition (ii) with the strict inequality $>-1$ instead of $\geq -1$. 
\item[(3)]
Suppose that $X^x_u=x<y=X^y_u$ and $0\leq X^x,X^y\leq 1$. Then, $X^y-X^x\geq 0$ and $\frac{1}{(X^y-X^x)_{-}}\ind_{\{(X^y-X^x)_->0\}}$ is $X^y-X^x$ integrable on $[u,\infty)$ if and only if the above condition (iii) is satisfied by $(\mathbf{F},\mathbf{Y})$ with $X=X^y, X'=X^x$. In addition, $X^y-X^x>0, (X^y-X^x)_->0$ on $[u,\infty)$ if and only if $(\mathbf{F},\mathbf{Y})$ and $X=X^y, X'=X^x$ satisfies the condition (iii) with the strict inequality $>-1$ instead of $\geq -1$.
\dbe
\el

\textbf{Proof.}
\ebe
\item[(1)]
Note that $d(1-Z)={^\bullet}(1-Z)d\widetilde{m}+dA$, also ${^\bullet}X^x=X^x_-$. From this we deduce
$$
\dcb
&&d((1-Z)_t-X^x_t)\\
&=&({^\bullet}(1-Z)_t-X^x_{t-})d\widetilde{m}_t-\mathbf{F}(X^x)_t^\top d\mathbf{Y}_t+dA_t\\

&=&({^\bullet}(1-Z)_t-X^x_{t-})\left(d\widetilde{m}_t-\frac{1}{{^\bullet}(1-Z)_t-X^x_{t-}}\ind_{\{{^\bullet}(1-Z)_t-X^x_{t-}\neq 0\}}\mathbf{F}(X^x)_t^\top d\mathbf{Y}_t\right)+dA_t\\

&=&{^\bullet}((1-Z)-X^x)_{t}\left(d\widetilde{m}_t-\frac{1}{{^\bullet}(1-Z)_t-X^x_{t-}}\ind_{\{{^\bullet}(1-Z)_t-X^x_{t-}\neq 0\}}\mathbf{F}(X^x)_t^\top d\mathbf{Y}_t\right)+dA_t.

\dce
$$
If $(\mathbf{F},\mathbf{Y})$ satisfies the condition (i), Lemma \ref{positive} is applicable, which yields the positivity of $(1-Z)-X^x$. Inversely, on the random interval $(S,T]$, by a direct computation, we have$$
\frac{\Delta_t(-M-X^x)}{{^\bullet}((1-Z)-X^x)_t}>-1,\ t\in(S,T],
$$
which is equivalent to say that the condition (i) is satisfied on $(S,T]$ by $(\mathbf{F},\mathbf{Y})$ with $X=X^x$.

\item[(2)]
The first assertion is the consequence of the Doleans-Dade exponential formula (cf. \cite[Theorem 9.39]{Yan}). The second assertion is the consequence of the same formula together with \cite[Theorem 2.62]{Yan}.

\item[(3)]
We make the same argument as in (2). \ok 
\dbe

\

As a consequence we obtain the following theorem.

\bethe\label{XZ}
Let $0<u<\infty$ and $0\leq x\leq 1-Z_u$. Then, the equation $(\natural_u)$ associated with a $\natural$-pair $(\mathbf{F},\mathbf{Y})$ has a unique solution $X^x$ on $[u,\infty)$ such that $0\leq X^x\leq 1-Z$. In particular, $X^x$ is a uniformly integrable $\mathbb{F}$ martingale on $[u,\infty)$. If $0\leq x\leq y\leq 1$, $X^x\leq X^y$ on $[u,\infty)$.
\ethe

\

The following theorem proves that the family of $\natural$-pair is not empty. 

\bethe\label{markov}
Let $m>0$ be an integer. Let $g(t,x)$ be any bounded continuously differentiable function defined on $\mathbb{R}_+\times \mathbb{R}$ taking values in $\mathbb{R}^m$. Let $\varphi$ be a $C^\infty$ increassing function on $\mathbb{R}_+$ such that $|\varphi(x)|\leq 2$ and $|\frac{\varphi(x)}{x}|\leq 1$. For $t\in \mathbb{R}_+$, we introduce the set $\mathsf{G}_{t}$ of $\mathbf{z}\in\mathbb{R}^m$ satisfying the two conditions:$$
\dcb
\circ&:&2\left|g(t,x)^\top\mathbf{z}\right|<1+\Delta_t\tilde{m},\mbox{ for $x\in\mathbb{R}$,}\\
\\
\circ\circ&:&\left[-\varphi'({^\bullet}(1-Z_{t-})-x)\varphi(x)g(t,x)+\varphi({^\bullet}(1-Z_{t-})-x)\varphi'(x)g(t,x)\right.\\
&&\hspace{3cm}+\left.\varphi({^\bullet}(1-Z_{t-})-x)\varphi(x)g'(t,x)\right]^\top\mathbf{z}
>-(1+\Delta_t\tilde{m}),\mbox{ for $x\in\mathbb{R}$}.
\dce
$$
(Here $g'(t,x)$ denotes the derivative with respect to $x$.) Then, for any $t\in\mathbb{R}_+$, the random set $\mathtt{G}_t$ is not empty, and the set-valued process $\mathtt{G}$ is $\mathbb{F}$ optional. There exists an $m$-dimensional $\mathbb{F}$ local martingale $\mathbf{Y}=(Y_1,\ldots,Y_m)$ whose jump at $t\in\mathbb{R}_+$, if it exists, is contained in $\mathsf{G}_{t}$. Let $$
\mathbf{F}(X)_t=f(t,X_{t-})=\varphi({^\bullet}(1-Z)_t-X_{t-})\varphi(X_{t-}) g(t,X_{t-}), \ X\in\mathbb{D}.
$$
Then, the above conditions (i), (ii) and (iii) with strict inequality $>-1$ instead of $\geq -1$ are satisfied for the pair $(\mathbf{F},\mathbf{Y})$.
\ethe

\textbf{Proof.} We recall that $\Delta_t\tilde{m}>-1$. Consequently, $\mathsf{G}_{t}$ contains always a no-empty neighbourhood of the origin. The optionality of $\mathsf{G}$ with respect to $\mathbb{F}$ can be proved in a usual way (cf. \cite{KK}). The existence of the $\mathbf{Y}$ is deduced from the measurable selection theorem. We now compute
$$
\dcb
&&\Delta_t \tilde{m}-\frac{1}{{^\bullet}(1-Z)_t-X_{t-}}\ind_{\{{^\bullet}(1-Z)_t-X_{t-}\neq 0\}}f(t,X_{t-})^\top\Delta_t\mathbf{Y}\\
&\geq&\Delta_t\tilde{m}-\ind_{\{{^\bullet}(1-Z)_t-X_{t-}\neq 0\}}2|g(t,X_{t-})^\top\Delta_t\mathbf{Y}|>-1,
\dce
$$
and
$$
\dcb
&&\Delta_t\tilde{m}+\frac{1}{X_{t-}}\ind_{\{X_{t-}\neq 0\}}f(t,X_{t-})^\top\Delta_t\mathbf{Y}\\
&\geq&\Delta_t\tilde{m}-\ind_{\{X_{t-}\neq 0\}}2|g(t,X_{t-})^\top\Delta_t\mathbf{Y}|>-1.\\

\dce
$$
These computations proves the conditions (i) and (ii) with strict inequality $>-1$. The condition (iii) with strict inequality is the sequence of the assumption $\circ\circ$, because then the map $$
x\longrightarrow x+x\Delta_t \tilde{m}+f(t,x)^\top\Delta_t\mathbf{Y}
$$
is strictly increasing.\ok

\subsection{The i$M_Z$ associated with the $\natural$-equation}

We assume always \textbf{Hy}$(Z)$.

\bethe Suppose $\mathcal{F}_\infty=\vee_{t\in\mathbb{R}_+}\mathcal{F}_t$. For $ 0<u<\infty$, let $(L^u_t:t\in[u,\infty))$ denote the solution of the
equation $(\natural_u)$ with the initial condition $L^u_u=1-Z_u$.
Set $L^u_\infty=\lim_{t\rightarrow\infty}L^u_t$. Set $$ \dcb
M^u_u=(1-Z_u),\\
M^u_t=\inf_{v\in\mathtt{Q},u<v\leq t}(L^v_t)^+\wedge(1-Z_t),\ t\in(u,\infty].\\
\dce
$$
Set finally$$
\dcb
M^0_t&=&\inf_{u\in\mathbb{Q}, 0<u\leq t}M^u_t,\ t\in(0,\infty],\\
M^0_0&=&\lim_{t\downarrow 0}M^0_t \ \mbox{ (which exists)},\\
M^\infty_t&=&1,\ \mbox{ for $t\in[0,\infty]$}.
\dce
$$
Then, for $0<u<\infty$, $M^u$ is $\mathbb{P}$ indistinguishable to $L^u$ on $[u,\infty]$, and $(M^u:0\leq u\leq \infty)$ is an i$M_Z$. 
\ethe 

The above i$\!M_Z$ will be said to be associated with the $\natural$-equation as well as the probability measure $\mathbb{Q}^\natural$ constructed in Theorem \ref{theoremIM} with this i$\!M_Z$ and with $\mathbb{P}$ will be said to be associated with the $\natural$-equation.

\textbf{Proof} As it has been proved in Theorem \ref{XZ}, $0\leq L^u_t\leq
(1-Z_t)$ for $0<u\leq t\leq \infty$, and $L^u_t\leq L^v_t$ for $0<
u<v\leq t\leq \infty$. The comparaison relation holds because $L^u,L^v$ satisfy the same equation ($\natural_v$) on $[v,\infty)$ and $L^u_v\leq (1-Z_v)=L^v_v$. Let $T$ be an $\mathbb{F}$ stopping time with $0< u\leq T\leq
\infty$. We have $$ 
\dcb
&&\mathbb{E}[M^u_T-L^u_T]\\
&=&\mathbb{E}[\mathbb{I}_{\{u<T\}}\inf_{v\in\mathtt{Q},u<v\leq T}L^v_T-\mathbb{I}_{\{u<T\}}L^u_T]\\
&=&\inf_{v\in\mathtt{Q},u<v}\mathbb{P}[\mathbb{I}_{\{u<T\}}L^{v\wedge T}_T]-\mathbb{E}[\mathbb{I}_{\{u<T\}}L^u_T]\\
&=&\inf_{v\in\mathtt{Q},u<v}\mathbb{P}[\mathbb{I}_{\{u<T\}}(1-Z_{v\wedge T})]-\mathbb{E}[\mathbb{I}_{\{u<T\}}(1-Z_u)]\\
&=&0.
\dce
$$
This shows that $M^u,L^u$ are indistinguishable on $[u,\infty]$
and in particular $M^u$ is càdlàg on $[u,\infty]$. 

Let $T$ be an $\mathbb{F}$ stopping time with $0< T\leq
\infty$. We have $$ 
\dcb
&&\mathbb{E}[M^0_T]\\
&=&\mathbb{E}[\inf_{v\in\mathtt{Q},0<v\leq T}L^v_T]\\
&=&\mathbb{E}[\inf_{v\in\mathtt{Q},v>0}L^{v\wedge T}_T]\\
&=&\inf_{v\in\mathtt{Q},v>0}\mathbb{E}[L^{v\wedge T}_T]\\
&=&\inf_{v\in\mathtt{Q},v>0}\mathbb{E}[(1-Z_{v\wedge T})]\\
&=&\mathbb{E}[(1-Z_{0})].
\dce
$$
The value $\mathbb{E}[M^0_T]$ does not depends on the stopping time $T$. According to \cite[Theorem 4.40]{Yan}, $M^0_t$ is a càdlàg $\mathbb{F}$ martingale on $(0,\infty]$. From this property we deduce that $M^0_0=\lim_{t\downarrow 0}M^0_t$ exists.

By definition, for
$0\leq t\leq \infty$, the map $u\in[0,t]\rightarrow M^u_t$ is a right
continuous non decreasing function. The theorem is proved. \ok

\

\section{Formula of enlargement of filtration in the case of Markovian coefficients}\label{efmc}

\subsection{The problem}

We consider the $\natural$-model associated with an $m$-dimensional $\natural$-pair $(\mathbf{F},\mathbf{Y})$. We suppose that the operator $\mathbf{F}(X)$ takes the particular form $\mathbf{F}(X)=f(\omega,t,X_{t-})$, where $f(\omega,t,x)$ is a map from $\Omega\times\mathbb{R}_+\times\mathbb{R}$ into $\mathbb{R}^m$ such that:
\ebe
\item
For $1\leq j\leq m$, for $\omega\in\Omega$ and $t\geq 0$, the partial derivative $\frac{\partial f_j}{\partial x}(\omega,t,x)$ exists and is bounded, which, for fixed $\omega$, is uniformly continuous with respect to $t$ in every compact set. 
\item
$f(\omega,t,0)=0$.
\dbe
We know by Theorem \ref{markov} that such a $\natural$-pair exists. It is customary to call such $\mathbf{F}$ a Markovian coefficient. 

Let $\mathbb{Q}^\natural$ be the probability measure on the product space associated the $\natural$-equation by Theorem \ref{theoremIM}. Let $\mathbb{G}=(\mathcal{G}_t)_{t\geq 0}$ be the filtration defined by $\cap_{s>t}(\mathcal{F}_{s}\vee \sigma (\tau \wedge
s))
$
completed with the $\mathbb{Q}^\natural$ negligible sets. In this section, we consider the problem of enlargement of filtration : whether the bounded $\mathbb{F}$ martingales are semimartingales in the filtration $\mathbb{G}$. 

In the computations below, the expectations are all taken under $\mathbb{Q}^\natural$ (recalling that $\mathbb{Q}^\natural$ coincides with $\mathbb{P}$ on $\mathcal{A}$).

\subsection{Preliminary results}

We need the following results.

\bl
For $0\leq u<\infty$, for $h$ a bounded Borel function on $[0,\infty]$, the process : $$
t\in[u,\infty]\rightarrow\int_0^u h(v) d_vM^v_t
$$ is a bounded $(\mathbb{P},\mathbb{F})$ martingale on $[u,\infty]$.  
\el

\textbf{Proof.} Consider the family of bounded Borel function $h$ on $[0,\infty]$ such that $$
\mathbb{E}[\int_0^u h(v) d_vM^v_T]
=\mathbb{E}[\int_0^u h(v) d_vM^v_u]
$$
for any $\mathbb{F}$ stopping times $T\geq u$. By monotone class theorem, we see that this family actually contains all bounded Borel function on $[0,\infty]$. According to \cite[Theorem 4.40]{Yan}, $\int_0^u h(v) d_vM^v_t, t\in[u,\infty)$ is a $(\mathbb{P},\mathbb{F})$ martingale. \ok

\bl\label{proj}
Let $V_s(\omega),0\leq s<\infty,\omega\in\Omega$, be a function such that, for fixed $s$, $V_s$ is $\mathcal{F}_\infty$ measurable and for fixed $\omega$, $s\rightarrow V_s(\omega)$ is càdlàg and non decreasing. We denote by $dV_s$ the induced measure on $[0,\infty)$. Let $F_s(t,\omega),0\leq s<\infty,0\leq t\leq \infty,\omega\in\Omega$, be a positive function measurable with respect to $\mathcal{B}[0,\infty]\otimes\mathcal{B}[0,\infty]\otimes\mathcal{F}_\infty$. Suppose $\mathbb{E}[\int_0^\infty dV_s]<\infty$. Then, $$
\mathbb{E}[\int_0^\infty F_sdV_s]
=\mathbb{E}[\int_0^\infty \left(\int_{[0,\infty]} F_s(v,\cdot)d_vM^v_\infty\right) dV_s].
$$
\el

\textbf{Proof.} By monotone class theorem, we need only to check
the relation for a function of form
$F_s(t,\omega)=h(s)H(t,\omega)$. We compute
\begin{eqnarray*}
\mathbb{E}[\int_0^\infty F_s dV_s]&=&
\mathbb{E}[H\int_0^\infty h(s)dV_s]\\
&=&\mathbb{E}[\int_{[0,\infty]} H(v,\cdot)d_vM^v_\infty\ \int_0^\infty h(s)dV_s]\\
&=&\mathbb{E}[\int_0^\infty \left(\int_{[0,\infty]}
F_s(v,\cdot)d_vM^v_\infty\right) dV_s]. 
\end{eqnarray*} The lemma
is proved. \ok

\bl\label{regular}
Suppose that there exists a prediction process $\chi$ of the conditional laws $\tau$ with respect to the filtration $\mathbb{F}$ under $\mathbb{Q}^\natural$ (cf. \cite[Theorem(13.1)]{aldous} and \cite{songthesis}). Then, for any $\mathbb{F}$ stopping time $0<T\leq \infty$, $\mathbb{Q}^\natural$ almost surely, $M^u_{T}=\chi_{T}[\tau\leq u]$ for all $u\in[0,T)$. On the other side, for fixed $u\in(0,\infty)$, for bounded function $h(\omega,s, v)$ measurable with respect to $\mathcal{P}(\mathbb{F})\times\mathcal{B}(\mathbb{R}_+)$, $$
\dcb
{^\bullet}\left(\int_{[0,u]} h(\cdot,v)\chi_{\infty}(dv)\right)
={^\bullet}\left(\int_{[0,u]} h(\cdot,v)d_vM^v_{\infty}\right)
=\int_{[0,u]} h(\cdot,v)\chi_{-}(dv)
\dce
$$
on $(u,\infty]$. In particular, $M^u_-=\chi_{-}[\tau\leq u]$ on $(u,\infty]$.
\el

\textbf{Proof.} We note that, $\mathbb{Q}^\natural$ almost surely, $$
M^u_{T}\ind_{\{u<T\}}=\mathbb{P}[\tau\leq u|\mathcal{F}_T]\ind_{\{u<T\}}=\chi_{T}[\tau\leq u]\ind_{\{u<T\}}
$$ 
for rational $u$. By the right continuity in $u$ of $M^u_{t}$ and of $\chi_{t}[\tau\leq u]$, we prove the first assertion. For the identity, let $\epsilon>0$ and $\varphi(x)$ be a non negative bounded continuous function on $\mathbb{R}_+$ whose support is contained in $[0,u+\epsilon]$. We consider $h(s,v)=H_sg(v)$, where $H$ is a bounded $\mathbb{F}$ predictable process and $g$ is a bounded continuous function on $\mathbb{R}_+$. We consider $T$ a $\mathbb{F}$ predictable stopping time such that $u+2\epsilon< T\leq \infty$ and $(T_n)_{n\geq 1}$ an non decreasing sequence of $\mathbb{F}$ stopping times such that $u+\epsilon< T_n<T$ and $T_n\uparrow T$. We compute
$$
\dcb
&&{^\bullet}\left(\int_{[0,u+\epsilon]} h(\cdot,v)\varphi(v)d_vM^v_{\infty}\right)_T\\
&=&\mathbb{E}[\int_{[0,u+\epsilon]} h(T,v)\varphi(v)d_vM^v_{\infty}|\mathcal{F}_{T-}]\\
&=&H_T\mathbb{E}[\int_{[0,u+\epsilon]} g(v)\varphi(v)d_vM^v_{\infty}|\mathcal{F}_{T-}]\\
&=&H_T\lim_{n\uparrow \infty}\mathbb{E}[\int_{[0,u+\epsilon]} g(v)\varphi(v)d_vM^v_{\infty}|\mathcal{F}_{T_n}]\\
&=&H_T\lim_{n\uparrow \infty}\int_{[0,u+\epsilon]} g(v)\varphi(v)d_vM^v_{T_n}\\
&=&H_T\lim_{n\uparrow \infty}\int_{[0,u+\epsilon]} g(v)\varphi(v)\chi_{T_n}(dv)\\
&=&\int_{[0,u+\epsilon]} h(T,v)\varphi(v)\chi_{T-}(dv).
\dce
$$ 
By the monotone class theorem, this identity remains valid for all bounded function $h(\omega,s, v)$ measurable with respect to $\mathcal{P}(\mathbb{F})\times\mathcal{B}(\mathbb{R}_+)$. Now, to finish the proof, let $\varphi$ decrease to $\ind_{[0,u]}$. \ok

\brem
\cite[Theorem(13.1)]{aldous} gives easy-to-use sufficient conditions to have a prediction process $\chi$. The reason to introduce $\chi$ is that $u\rightarrow \chi_{t-}[\tau\leq u]$ is a true distribution function, whilst the behavior of $u\rightarrow M^u_{t-}$ is unknown.
\erem

\bl\label{develop}
Let $0\leq s,u<\infty$. Let $A\in\mathcal{F}_s$. Let $X$ be a bounded $\mathbb{F}$ martingale such that $X_{s\vee u}=0$ and $\int_{s\vee u}^\infty |d\cro{\tilde{m},X}| 
+ \int_{s\vee u}^\infty |d\cro{\mathbf{Y},X}|$ is integrable. Let $\mathtt{p}(\omega,w,v)$, $(\omega,w,v)\in\Omega\times\mathbb{R}_+\times\mathbb{R}_+$, denote the function $$
\mathtt{p}(\omega,w,v)=\left(\frac{\partial f}{\partial x}(w,\chi_{w-}[\tau\leq v])\ind_{\{\chi_{w-}[\tau=v]=0\}}
+\frac{f(w,\chi_{w-}[\tau\leq v])-f(w,\chi_{w-}[\tau<v])}{\chi_{w-}[\tau=v]}\ind_{\{\chi_{w-}[\tau=v]>0\}}\right).
$$
We have$$
\dcb
\mathbb{E}[\ind_A \ind_{\{\tau\leq u\}} X_\infty]
&=&\mathbb{E}[\ind_A \ind_{\{\tau\leq u\}}\left(\int_{s\vee u}^\infty d\cro{\tilde{m},X}_w 
+ \int_{s\vee u}^\infty \mathtt{p}(w,\tau)^\top d\cro{\mathbf{Y},X}_w\right)].
\dce
$$
\el

\textbf{Proof.} Let us denote the $\mathbb{F}$ predictable bracket of two $\mathbb{F}$ local martingales $X,X'$ by $\cro{X,X'}$ (when it exists). Suppose firstly $u>0$. We compute$$
\dcb
\mathbb{E}[\ind_A \ind_{\{\tau\leq u\}} X_\infty]
&=&\mathbb{E}[\ind_A M^u_\infty X_\infty]\\
&=&\mathbb{E}[\ind_A \int_{s\vee u}^\infty d\cro{M^u,X}_w]\\

&=&\mathbb{E}[\ind_A \left(\int_{s\vee u}^\infty M^u_{w-} d\cro{\tilde{m},X}_w 
+ \int_{s\vee u}^\infty f(w,M^u_{w-})^\top d\cro{\mathbf{Y},X}_w\right)].
\dce
$$
By Lemma \ref{regular}, we can write, for all $w>s\vee u$,
$$
\dcb
&&f(w,M^u_{w-})\\
&=&f(w,M^0_{w-})+f(w,M^u_{w-})-f(w,M^0_{w-})\\
&=&f(w,\chi_{w-}[\tau\leq 0])+f(w,\chi_{w-}[\tau\leq u])-f(w,\chi_{w-}[\tau\leq 0])\\
&=&\frac{f(w,\chi_{w-}[\tau= 0])}{\chi_{w-}[\tau= 0]}\ind_{\{\chi_{w-}[\tau= 0]>0\}}+\int_{(0,u]}\frac{\partial f}{\partial x}(w,\chi_{w-}[\tau<v])\chi_{w-}(dv)\\
&&+\sum_{0<v\leq u}(f(w,\chi_{w-}[\tau\leq v])-f(w,\chi_{w-}[\tau<v])-\frac{\partial f}{\partial x}(w,\chi_{w-}[\tau<v])\chi_{w-}[\tau=v])\\

&=&\int_{[0,u]}\left(\frac{\partial f}{\partial x}(w,\chi_{w-}[\tau\leq v])\ind_{\{\chi_{w-}[\tau=v]=0\}}
+\frac{f(w,\chi_{w-}[\tau\leq v])-f(w,\chi_{w-}[\tau<v])}{\chi_{w-}[\tau=v]}\ind_{\{\chi_{w-}[\tau=v]>0\}}\right) \chi_{w-}(dv)\\

&=&\int_{[0,u]} \mathtt{p}(w,v) \chi_{w-}(dv).
\dce
$$
We note that $\mathtt{p}$ is bounded and is $\mathcal{P}(\mathbb{F})\times\mathcal{B}(\mathbb{R}_+)$ measurable. Applying Lemma \ref{proj} and Lemma \ref{regular} we continue the computation:
$$
\dcb
\mathbb{E}[\ind_A \ind_{\{\tau\leq u\}} X_\infty]

&=&\mathbb{E}[\ind_A \left(\int_{s\vee u}^\infty M^u_{w-} d\cro{\tilde{m},X}_w + \int_{s\vee u}^\infty \int_{[0,u]} \mathtt{p}(w,v)^\top \chi_{w-}(dv)d\cro{\mathbf{Y},X}_w\right)]\\
\\

&=&\mathbb{E}[\ind_A \left(\int_{s\vee u}^\infty M^u_{\infty} d\cro{\tilde{m},X}_w 
+ \int_{s\vee u}^\infty \int_{[0,u]} \mathtt{p}(w,v)^\top \chi_{\infty}(dv)d\cro{\mathbf{Y},X}_w\right)]\\
\\

&=&\mathbb{E}[\ind_A \left(\int_{s\vee u}^\infty \ind_{\{\tau\leq u\}} d\cro{\tilde{m},X}_w 
+ \int_{s\vee u}^\infty \int_{[0,u]} \mathtt{p}(w,v)^\top d_vM_{\infty}^vd\cro{\mathbf{Y},X}_w\right)]\\
\\

&=&\mathbb{E}[\ind_A \left(\int_{s\vee u}^\infty \ind_{\{\tau\leq u\}} d\cro{\tilde{m},X}_w 
+ \int_{s\vee u}^\infty \ind_{\{\tau\leq u\}}\mathtt{p}(w,\tau)^\top d\cro{\mathbf{Y},X}_w\right)]
\dce
$$
Now let $u\downarrow 0$ to finish the proof of the lemma. \ok

\subsection{The proof of the enlargement of filtration formula}

Recall that the canonical decomposition of $Z$ is given by $Z=M-A$. Recall the result obtained in \cite{J,JY2}.

\bl\label{Jlemma} Let $X$ be a
bounded $(\mathbb{P},\mathbb{F})$ martingale and $B^X$ the
$(\mathbb{Q}^\natural,\mathbb{F})$ predictable dual projection of the jump
process $t\rightarrow \Delta X_\tau\ind_{\{0<\tau\leq t\}}$. Then, $$
X_{\cdot\wedge \tau} - \int_0^{\cdot\wedge
\tau}\frac{1}{Z_{s-}}(d\cro{M,X}_s+dB^X_s)
$$ is a $(\mathbb{Q}^\natural,\mathbb{G})$  local martingale.
\el

\bethe \label{decompositionFormula}
Let $X$ be a bounded $(\mathbb{P},\mathbb{F})$ martingale. Then the process
$$
\dcb
X_{t}-X_0-\int_{0}^{t}\ind_{\{s\leq \tau\}}\frac{1}{Z_{s-}}(d\cro{M,X}_s+dB^X_s)\\
\\
\hspace{1.5cm}+\int_{0}^{u} \ind_{\{\tau<w\}}\frac{1}{1-{^\bullet}Z_w}d\cro{M,X}_w 
- \int_{0}^{u} \ind_{\{\tau<w\}}\mathtt{p}(w,\tau)^\top d\cro{\mathbf{Y},X}_w\\
\dce
$$ is a $(\mathbb{Q}^\natural,\mathbb{G})$ local martingale. \ethe

\textbf{Proof.} We write $$
X_t-X_0=X_{\tau\vee t}-X_\tau+X_{\tau\wedge t}-X_0.
$$
The $\mathbb{G}$ semimartingale decomposition for $X_{\tau\wedge t}-X_0$ is given in Lemma \ref{Jlemma}. We need only to prove the formula for $X_{\tau\vee t}-X_\tau$. Without loss of generality we suppose that $X$ is stopped so that everything in the computations below is integrable. Let $0< s<t<\infty,0< u<\infty$ and $A\in\mathcal{F}_s$.
$$
\dcb
&&\mathbb{E}[\ind_A\ind_{\{\tau\leq u\}}(X_{\tau\vee t}-X_{\tau\vee s})]\\
&=&\mathbb{E}[\ind_A\ind_{\{\tau=0\}}(X_{t}-X_{s})]+\lim_{n\rightarrow\infty}\mathbb{E}[\ind_A\sum_{k=1}^n\ind_{\{\frac{(k-1)u}{n}<\tau\leq \frac{ku}{n}\}}(X_{\tau\vee t}-X_{\tau\vee s})]\\
&=&\mathbb{E}[\ind_A\ind_{\{\tau=0\}}(X_{t}-X_{s})]+\lim_{n\rightarrow\infty}\mathbb{E}[\ind_A\sum_{k=1}^n\ind_{\{\frac{(k-1)u}{n}<\tau\leq \frac{ku}{n}\}}(X_{\frac{ku}{n}\vee t}-X_{\frac{ku}{n}\vee s})].
\dce
$$
Set $t_n=\frac{ku}{n}\vee
t$ and $s_n=\frac{ku}{n}\vee s$. By Lemma \ref{develop}, we have$$
\dcb
&&\mathbb{E}[\ind_A\ind_{\{\frac{(k-1)u}{n}<\tau\leq \frac{ku}{n}\}}(X_{\frac{ku}{n}\vee t}-X_{\frac{ku}{n}\vee s})]\\

&=&\mathbb{E}[\ind_A\ind_{\{\tau\leq \frac{ku}{n}\}}(X_{\frac{ku}{n}\vee t}-X_{\frac{ku}{n}\vee s})]
-\mathbb{E}[\ind_A\ind_{\{\tau\leq \frac{(k-1)u}{n}\}}(X_{\frac{ku}{n}\vee t}-X_{\frac{ku}{n}\vee s})]\\

&=&\mathbb{E}[\ind_A \ind_{\{\tau\leq \frac{ku}{n}\}}\left(\int_{s_n}^{t_n} d\cro{\tilde{m},X}_w 
+ \int_{s_n}^{t_n} \mathtt{p}(w,\tau)^\top d\cro{\mathbf{Y},X}_w\right)]\\
&&-\mathbb{E}[\ind_A \ind_{\{\tau\leq \frac{(k-1)u}{n}\}}\left(\int_{s_n}^{t_n} d\cro{\tilde{m},X}_w 
+ \int_{s_n}^{t_n} \mathtt{p}(w,\tau)^\top d\cro{\mathbf{Y},X}_w\right)]\\

&=&\mathbb{E}[\ind_A \ind_{\{\frac{(k-1)u}{n}<\tau\leq \frac{ku}{n}\}}\left(\int_{s_n}^{t_n} d\cro{\tilde{m},X}_w 
+ \int_{s_n}^{t_n} \mathtt{p}(w,\tau)^\top d\cro{\mathbf{Y},X}_w\right)].
\dce
$$
Also,
$$
\dcb
\mathbb{E}[\ind_A\ind_{\{\tau=0\}}(X_{t}-X_{s})]
=\mathbb{E}[\ind_A\ind_{\{\tau=0\}}\left(\int_{s}^{t} d\cro{\tilde{m},X}_w 
+ \int_{s}^{t} \mathtt{p}(w,\tau)^\top d\cro{\mathbf{Y},X}_w\right)]
\dce
$$
We turn back to $\mathbb{E}[\ind_A\ind_{\{\tau\leq u\}}(X_{\tau\vee t}-X_{\tau\vee s})]$:
$$
\dcb
&&\mathbb{E}[\ind_A\ind_{\{\tau\leq u\}}(X_{\tau\vee t}-X_{\tau\vee s})]\\
&=&\mathbb{E}[\ind_A\ind_{\{\tau=0\}}(X_{t}-X_{s})]+\lim_{n\rightarrow\infty}\mathbb{E}[\ind_A\sum_{k=1}^n\ind_{\{\frac{(k-1)u}{n}<\tau\leq \frac{ku}{n}\}}(X_{\frac{ku}{n}\vee t}-X_{\frac{ku}{n}\vee s})]\\

&=&\mathbb{E}[\ind_A\ind_{\{\tau=0\}}\left(\int_{s}^{t} d\cro{\tilde{m},X}_w 
+ \int_{s}^{t} \mathtt{p}(w,\tau)^\top d\cro{\mathbf{Y},X}_w\right)]\\

&&+\lim_{n\rightarrow\infty}\mathbb{E}[\ind_A\sum_{k=1}^n\ind_{\{\frac{(k-1)u}{n}<\tau\leq \frac{ku}{n}\}}\left(\int_{s_n}^{t_n} d\cro{\tilde{m},X}_w 
+ \int_{s_n}^{t_n} \mathtt{p}(w,\tau)^\top d\cro{\mathbf{Y},X}_w\right)]\\

&=&\mathbb{E}[\ind_A\ind_{\{\tau=0\}}\left(\int_{s}^{t} d\cro{\tilde{m},X}_w 
+ \int_{s}^{t} \mathtt{p}(w,\tau)^\top d\cro{\mathbf{Y},X}_w\right)]\\

&&+\mathbb{E}[\ind_A\ind_{\{0<\tau\leq u\}}\left(\int_{\tau\vee s}^{\tau\vee t} d\cro{\tilde{m},X}_w 
+ \int_{\tau\vee s}^{\tau\vee t} \mathtt{p}(w,\tau)^\top d\cro{\mathbf{Y},X}_w\right)]\\

&=&\mathbb{E}[\ind_A\ind_{\{\tau\leq u\}}\left(\int_{s}^{t}\ind_{\{\tau<w\}} d\cro{\tilde{m},X}_w 
+ \int_{s}^{t} \ind_{\{\tau<w\}}\mathtt{p}(w,\tau)^\top d\cro{\mathbf{Y},X}_w\right)].
\dce
$$
Note that, for $0<s<s'$, $\mathcal{G}_s\subset \mathcal{F}_{s'}\otimes\sigma(\tau)$. The above formula implies$$
\mathbb{E}[\ind_B(X_{\tau\vee t}-X_{\tau\vee s'})]
=\mathbb{E}[\ind_B\left(\int_{s'}^{t} \ind_{\{\tau<w\}}d\cro{\tilde{m},X}_w 
+ \int_{s'}^{t} \ind_{\{\tau<w\}}\mathtt{p}(w,\tau)^\top d\cro{\mathbf{Y},X}_w\right)]
$$
for any $B\in \mathcal{G}_s$. Taking the limit when $s'\downarrow s$, we prove the theorem. \ok

\

\section{Regularity of $u\rightarrow M^u$ in the case of Markovian coefficients}\label{conditionallaw}

Consider the $\natural$-equation associated with a $\natural$-pair $(\mathbf{F},\mathbf{Y})$ of the type in Theorem \ref{markov}:$$
\dcb
\Delta_tY&\in&\mathtt{G}_t,\\
\mathbf{F}(X)_t&=&\varphi({^\bullet}(1-Z)_t-X_{t-})\varphi(X_{t-}) g(t,X_{t-}), \ X\in\mathbb{D}.
\dce
$$ 
We suppose moreover that $\varphi(x)=x$, for $x\in[0,1]$; the function $g$ is autonomuous, i.e. $g(t,x)=g(x)$ (in the sense of \cite{protter}); $g$ is $C^\infty$ with bounded derivatives of all order; $g$ has a compact support. Consider the i$\!M_Z=(M^u:0\leq u\leq \infty)$ associated with $(\mathbf{F},\mathbf{Y})$. Since $0\leq M^u\leq 1-Z$, we have$$
\varphi({^\bullet}(1-Z)-M^u_-)\varphi(M^u_-)=({^\bullet}(1-Z)-M^u_-)M^u_-.
$$
We consider then the stochastic differential equation $$
\left\{\dcb
dX_t=X_{t-}d\tilde{m}_t+({^\bullet}(1-Z)_t-X_{t-})X_{t-}g(X_{t-})^\top d\mathbf{Y}_t,\ u\leq t<\infty,\\
X_u=x.
\dce
\right.
$$
Theorem 39 and Theorem 65 in \cite[Chapter V Section 10]{protter} are applicable to such an equation. Let $x\rightarrow \Xi_t^{u}(x)$ be the stochastic differential flow defined by this equation. We note that, for $0<u<\infty$, $M^u$ satisfies this equation (for $M^u$ the explosion time being $\infty$).

\bethe\label{regularity}
Let $0< t<\infty$. Let$$
\kappa_v=(1+\Delta_v\widetilde{m}-(1-Z_{v-})g(1-Z_{v-})^\top\Delta_v\mathbf{Y}), \ \mbox{ for $v\in(0,t]$}.
$$
Then, when $\Delta_vA>0, v\in(0,t]$, $$
M^v_t-M^{v-}_t=\Xi_t^v(1-Z_v)-\Xi_t^v(1-Z_v-\kappa_v\Delta_vA).
$$
When $\Delta_vA=0$,$$
\dcb
\lim_{u\uparrow v}\frac{M^v_t-M^u_t}{A_v-A_u}&=\frac{d\Xi_t^v}{dx}(1-Z_v)\kappa_v,\ &\mbox{ for $v\in(0,t]$},\\
\\
\lim_{u\downarrow v}\frac{M^v_t-M^u_t}{A_v-A_u}&=\frac{d\Xi_t^v}{dx}(1-Z_v),\ &\mbox{ for $v\in(0,t)$}.
\dce
$$
\ethe

\textbf{Proof.} For any $0< t<\infty$, we consider the difference $M^v_t-M^u_t$ for $u,v\in(0,t], u<v$. By the uniqueness of stochastic differential equation,
$$
\dcb
M^v_t-M^u_t
&=&\Xi_t^v(1-Z_v)-\Xi_t^v(M^u_v)
=\frac{d\Xi_t^v}{dx}(\xi)((1-Z_v)-M^u_v),
\dce
$$
where $\xi$ is some random variable such that $M^u_v\leq \xi\leq 1-Z_v$. The process $(1-Z_v)-M^u_v$ satisfies the equation $$
\left\{
\dcb
d((1-Z)-M^u)_t&=&{^\bullet}((1-Z)-M^u)_tdW^u_t+dA_t,\ t\in[u,\infty),\\
\\
((1-Z)-M^u)_u&=&0,
\dce
\right.
$$
with$$
\dcb
dW^u_t
&=&d\tilde{m}_t-M^u_{t-}g(M^u_{t-})^\top d\mathbf{Y}_t.
\dce
$$
According to Lemma \ref{separable}, $$
(1-Z_t)-M^u_t=\mathcal{E}(\ind_{(u,\infty)}\centerdot W^u)_t\int_u^t\frac{1}{\mathcal{E}(\ind_{(u,\infty)}\centerdot W^u)_{s-}}dA_s,\ t\in[u,\infty).
$$
By the property of the i$\!M_Z$, $
\lim_{u\uparrow t}M^u_{t-}=1-Z_{t-}.
$
By Doleans-Dade formula, $$
\lim_{u\uparrow t}\mathcal{E}(\ind_{(u,\infty)}\centerdot W^u)_{t-}=1.
$$ 
Consequently,$$
\dcb
\lim_{u\uparrow t}\mathcal{E}(\ind_{(u,\infty)}\centerdot W^u)_{t}
&=&\lim_{u\uparrow t}\mathcal{E}(\ind_{(u,\infty)}\centerdot W^u)_{t-}(1+\Delta_t\widetilde{m}-M^u_{t-}g(M^u_{t-})^\top\Delta_t\mathbf{Y})\\
\\
&=&(1+\Delta_t\widetilde{m}-(1-Z_{t-})g(1-Z_{t-})^\top\Delta_t\mathbf{Y}) = \kappa_t.
\dce
$$
We also need to estimate the difference, for $0<u<v\leq t$,$$
\dcb
\frac{d\Xi_t^v}{dx}(y)-\frac{d\Xi_t^u}{dx}(y)
&=&\frac{d\Xi_t^v}{dx}(y)-\frac{d\Xi_t^v}{dx}(\Xi_v^u(y))\frac{d\Xi_v^u}{dx}(y)\\
\\
&=&\frac{d\Xi_t^v}{dx}(y)-\frac{d\Xi_t^v}{dx}(\Xi_v^u(y))+\frac{d\Xi_t^v}{dx}(\Xi_v^u(y))(1-\frac{d\Xi_v^u}{dx}(y)).
\dce
$$
When $u$ is fixed and $v\downarrow u$, by \cite[Chapter V Section 7 Theorem 39]{protter}, we have $\lim_{v\downarrow u}\Xi_v^u(y)=y$ and $\lim_{v\downarrow u}\frac{d\Xi_v^u}{dx}(y)=1$. According to the same reference, by Kolmogorov continuity criterion, for any $a\in(0,t)$, the random map $y\rightarrow \frac{d\Xi_t^v}{dx}(y)$ is continuous uniformly with respect to $v\in[a,t]$. This uniform continuity implies$$
\lim_{v\downarrow u}\frac{d\Xi_t^v}{dx}(y)-\frac{d\Xi_t^u}{dx}(y)=0.
$$
Now we write 
$$
\dcb
\frac{M^v_t-M^u_t}{A_v-A_u}
&=&\frac{d\Xi_t^v}{dx}(\xi)((1-Z_v)-M^u_v)\frac{1}{A_v-A_u}\\
&=&\frac{d\Xi_t^v}{dx}(1-Z_u)\frac{(1-Z_v)-M^u_v}{A_v-A_u}+(\frac{d\Xi_t^v}{dx}(\xi)-\frac{d\Xi_t^v}{dx}(1-Z_u))\frac{(1-Z_v)-M^u_v}{A_v-A_u}\\

&=&\frac{d\Xi_t^v}{dx}(1-Z_v)\frac{(1-Z_v)-M^u_v}{A_v-A_u}+(\frac{d\Xi_t^v}{dx}(\xi)-\frac{d\Xi_t^v}{dx}(1-Z_v))\frac{(1-Z_v)-M^u_v}{A_v-A_u}

\dce
$$
We consider separately limits when $v$ is fixed and $u\uparrow v$ or when $u$ is fixed and $v\downarrow u$, we have$$
\dcb
\lim_{u\uparrow v}((1-Z_v)-M^u_v)&=\kappa_v\Delta_vA,&\\

\lim_{u\uparrow v}\frac{(1-Z_v)-M^u_v}{A_v-A_u}&= \kappa_v,\\

\lim_{v\downarrow u}((1-Z_v)-M^u_v)&=0,\\

\lim_{v\downarrow u}\frac{(1-Z_v)-M^u_v}{A_v-A_u}&= 1.\\
\dce
$$
From these computations we deduce that$$
\dcb
\lim_{u\uparrow v}\frac{M^v_t-M^u_t}{A_v-A_u}&=\frac{\Xi_t^v(1-Z_v)-\Xi_t^v(1-Z_v-\kappa_v\Delta_vA)}{\Delta_vA},& \mbox{ if $\Delta_vA>0$, $v\in(0,t]$,}\\

\lim_{v\downarrow u}\frac{M^v_t-M^u_t}{A_v-A_u}&=0,& \mbox{ if $\Delta_uA>0$,  $u\in(0,t)$,}\\

\lim_{u\uparrow v}\frac{M^v_t-M^u_t}{A_v-A_u}&=\frac{d\Xi_t^v}{dx}(1-Z_v)\kappa_v,& \mbox{ if $\Delta_vA=0$, $v\in(0,t]$,}\\

\lim_{v\downarrow u}\frac{M^v_t-M^u_t}{A_v-A_u}&=\frac{d\Xi_t^u}{dx}(1-Z_u),& \mbox{ if $\Delta_uA=0$, $u\in(0,t)$.} \ok
\dce
$$

\

The knowledge on the regularity of the random map $u\rightarrow M^u_t$ for $u\in[0,t]$ will be very helpful in the study of the non-arbitrage property or of the optimization problem. In the immediate term, this regularity already gives answers to several questions about the $\natural$-model.

In \cite{JS2} we had proved results under the assumption that $u\in(0,t]\rightarrow M^u_t$ was continuous. We asked if this continuity effectively hold in the continuous $\natural$-model. Now we can say:

\bcor
For the $\natural$-model in Theorem \ref{regularity}, for any $0<t<\infty$, the map $u\in(0,t]\rightarrow M^u_t$ is continuous if and only if $A$ is continuous.
\ecor

Another question was if $d_uM^u_t$ is linked with $dA_u$. If we refine slightly the proof of Theorem \ref{regularity} we see that, for any $0<a<b<\infty$, for almost all $\omega$, the quotient $\frac{M^v_t(\omega)-M^u_t(\omega)}{A_v(\omega)-A_u(\omega)}$ is uniformly bounded for $v,u\in[a,b]$. This implies:

\bcor
For the $\natural$-model in Theorem \ref{regularity}, for any $0<t<\infty$, $d_uM^u_t$ is absolutely continuous with respect to $dA_u$ on $(0,t]$.
\ecor

The absolute continuity property of $d_uM^u_t$ with respect to $dA_u$ is established for finite interval $(0,t], t<\infty$. The situation when $t=\infty$ can be very different. We have the following lemma from \cite{aksamit}:

\bl
Suppose that $Z_t=e^{-t}, t\in[0,\infty)$. Let $B$ be a linear Brownian motion. Consider the $\natural$-equation:$$
(\natural_u):\left\{
\dcb
dX_t&=&X_t(1-e^{-t}-X_t)dB_t,\ t\in[u,\infty),\\
\\
X_u&=&1-e^{-u},
\dce
\right.
$$
for $0<u<\infty$, and the associated i$\!M_Z=(M^u:0\leq u\leq \infty)$. Then, the map $u\in[0,\infty]\rightarrow M^u_\infty$ takes only two values $0$ or $1$.
\el

In this lemma, clearly $d_uM^u_\infty$ is not absolutely continuous with respect to $dA_u=e^{-u}du$.

To end this section, notice that, if a one-default model satisfies the density hypothesis with respect to a determinist measure $\mu$ in the sense of \cite{ejj}, the random map $u\rightarrow M^u_t$ for $u\in[0,t]$ is absolutely continuous with respect to $\mu$. Consequently, we have: 

\bcor
There exists a $\natural$-model which does not satisfies the density hypothesis.
\ecor

\

\


\begin{thebibliography}{99}

\bibitem{aksamit} Aksamit A. \textit{Rôle de l'information dans les marchés financiers} Thesis University Evry (2013) 

\bibitem{ACDJ} Aksamit A. and Choulli T. and Deng J. and Jeanblanc M. "Non-arbitrage up to random horizon and after honest times for semimartingale models" working paper (2013)

\bibitem{aldous} Aldous D. \textit{Weak convergence and the general theory of processes} University of California (1981)


\bibitem{lisbonn} Bielecki T.R. and Jeanblanc M. and Rutkowski M. \textit{Credit Risk} Lisbonn Lecture (2006)

\bibitem{BJR} Bielecki T.R. and Jeanblanc M. and Rutkowski M.  \textit{Credit Rist Modelling} Osaka University Press (2009)



\bibitem{CJZ} Callegaro G. and Jeanblanc M. and Zargari B. "Carthaginian enlargement of filtrations" \textit{ESAIM: Probability and Statistics} \textbf{17} 550-566 (2013)


\bibitem{ejj} El Karoui N. and Jeanblanc M. and Ying J.  "What happens after a default : the conditional density approach" \textit{Stochastic Processes and their Applications} \textbf{120} 1011-1032 (2010)

\bibitem{EJY} Elliot R.J. and Jeanblanc M. and Yor M. "On models of default risk" \textit{Mathematical Finance} \textbf{10} 179-195 (2000)

\bibitem{FJS} Fontana C. and Jeanblanc M. and Song S. "On arbitrages arising from honest times" arxiv:1207.1759 (2013)



\bibitem{Yan} He S.W., Wang J.G., Yan J.A. \textit{Semimartingale Theory And Stochastic Calculues} Science Press CRC Press Inc (1992)


\bibitem{jacod} Jacod J. \textit{Calcul Stochastique et Problèmes de Martingales} Lecture Notes in Mathematics \textbf{714} Springer (1979)

\bibitem{JLC} Jeanblanc M. and LeCam Y. "Progressive enlargement of filtrations with initial times" \textit{Stochastic Processes and their Applications} \textbf{119} 2523-2543 (2009)

\bibitem{JS1}Jeanblanc M. and Song S. "An explicit  model of default time with given survival probability" \textit{Stochastic Processes and their Applications} \textbf{121}(8) 1678-1704 (2010)

\bibitem{JS2}Jeanblanc M. and Song S. "Random times with given survival probability and their $\mathbb{F}$-martingale decomposition formula" \textit{Stochastic Processes and their Applications} \textbf{121}(6) 1389-1410 (2010)


\bibitem{J} Jeulin T. \textit{Semi-martingales et grossissement d'une filtration} Lecture Notes in Mathematics \textbf{833} Springer (1980)

\bibitem{JY}Jeulin T. and Yor M. "Nouveaux résultats sur le grossissement des tribus". \textit{Ann. Scient. Ec. Norm. Sup.} \textbf{4} t.11 429-443 (1978)

\bibitem{JY2}Jeulin T. and Yor M. "Grossissement d'une filtration et semi-martingales: formules explicites" \textit{Séminaire de Probabilités} \textbf{12} 78-97 (1978)

\bibitem{KK} Karatzas I. and Kardaras C. "The numéraire portfolio in semimartingale financial models" \textit{Finance and Stochastics} \textbf{11}(4) 447-493 (2007)


\bibitem{knight} Knight F. \textit{Essays on the prediction process} IMS Lecture Notes Monograph Series (1981)



\bibitem{LR} Li L. and Rutkowski M. "Random times and multiplicative systems" \textit{Stochastic Processes and their Applications} \textbf{122}(5) 2053-2077 (2012)

\bibitem{MeYo} Meyer P. and Yoeurp C. "Sur la décomposition multiplicative des sousmartingales positives" \textit{Séminaire de Probabilités} \textbf{10} 501-504 (1976)


\bibitem{MYr} Meyer P. and Yor M. "Sur la théorie de la prédiction" \textit{Séminaire de Probabilités} \textbf{10} 104-117 (1976)

\bibitem{NP} Nikeghbali A. and Platen E. "On honest times in financial modeling" arxiv:0808.2892 (2008)

\bibitem{NY} Nikeghbali A. and Yor M. "A definition and some characteristic properties of pseudo-stopping times" \textit{The Annals of Probability} \textbf{33}(5) 1804-1824 (2005)

\bibitem{protter}
\textsc{Protter P.} \textit{Stochastic integration and differential equations} Springer (2005)

\bibitem{RY} Revuz D. and Yor M. \textit{Continuous martingales and Brownian motion} Springer (2004)

\bibitem{songthesis} Song S. \textit{Grossissement d'une filtration et problèmes connexes} Thesis Université Paris VI (1987)

\bibitem{songdrift} Song S. "Drift operator in a market affected by the expansion of information flow: a case study" arxiv:1207.1662 (2012)



\bibitem{yoeurp} Yoeurp C. "Décomposition des martingales locales et formules exponentielles" \textit{Séminaire de Probabilités} \textbf{10} 432-480 (1976)

\bibitem{yor} Yor M. "Sur les théories du filtrage et de la prédiction" \textit{Séminaire de Probabilités} \textbf{11} 257-297 (1977)


\end{thebibliography}
\end{document}